\documentclass{amsart}
\usepackage{amsthm,amssymb, amsmath, amscd}
\usepackage{tikz, tikz-cd}
\usepackage{graphicx}  
\usepackage{verbatim}
\usepackage{hyperref}
\usepackage[a4paper]{geometry}

\newtheorem{theorem}{Theorem}

\theoremstyle{definition}

\makeatletter
\def\x@arrow{\DOTSB\Relbar}
\def\xlongequalsignfill@{\arrowfill@\x@arrow\Relbar\x@arrow}
\newcommand{\xlongequal}[2][]{%
        \ext@arrow 0099\xlongequalsignfill@{#1}{#2}}
\def\xlongrightarrowfill@{\arrowfill@\relbar\relbar\longrightarrow}
\newcommand{\xlongrightarrow}[2][]{%
        \ext@arrow 0099\xlongrightarrowfill@{#1}{#2}}
        
\begin{document}

\title{Cohomologies on hypercomplex manifolds}
\author{Mehdi Lejmi}
\address{Department of Mathematics, Bronx Community College of CUNY, Bronx, NY 10453, USA.}
\email{mehdi.lejmi@bcc.cuny.edu}
\author{Patrick Weber}
\address{D\'epartement de Math\'ematiques, Universit\'e libre de Bruxelles CP218, Boulevard du Triomphe, Bruxelles 1050, Belgique. }
\email{pweber@ulb.ac.be}

%%%%%%%%%%%%%%%%%%%%%%%%%%%%%%%%%%%%%%%%%%%%%%%%%%%%%%%%%%%%

\begin{abstract}
We review some cohomological aspects of complex and hypercomplex manifolds and underline the differences between both realms. Furthermore, we try to highlight the similarities between compact complex surfaces on one hand and compact hypercomplex manifolds of real dimension~8 with holonomy of the Obata connection in~$\mathrm{SL}(2,\mathbb{H})$ on the other hand.
\end{abstract}

\maketitle

%%%%%%%%%%%%%%%%%%%%%%%%%%%%%%%%%%%%%%%%%%%%%%%%%%%%%%%%%%%%

\section{Introduction}
\label{sec:1}

We describe a recipe that allows one to adapt some cohomological results from complex manifolds to hypercomplex manifolds. A hypercomplex manifold is a complex manifold together with a second complex structure that anticommutes with the first one. To extract cohomological information out of a hypercomplex manifold, we may thus start with the double complex of the underlying complex manifold, twist this data by the second complex structure and see what information we get about the hypercomplex manifold in question. This approach turns out to be surprisingly successful if we want to adapt results from complex geometry to hypercomplex geometry and the resulting cohomology groups have the additional advantage of being easily computable.

\hfill

We would like to anticipate that this way of proceeding also suffers from some drawbacks and that there is an alternative approach available in the literature. If a manifold admits two anticommuting complex structures $I$ and~$J$, then $K=IJ$ is another almost-complex structure, anticommuting with both $I$ and~$J$. This then leads to a whole 2-sphere worth of almost-complex structures 
$$S^2= \{a I + bJ + cK~|~a^2+b^2+c^2 =1\}$$
and it has been shown that all of these almost-complex structures are integrable as soon as $I$ and $J$ are (see for example~\cite{kaledin}). From this point of view, all these complex structures should be treated equally on a hypercomplex manifold and singling out a preferred complex structure, as we do with~$I$, is not very natural. Unfortunately, the cohomology groups based upon the ``averaged complex structures" often tend to be quite cumbersome to work with and less explicit to compute. For further information, we refer the interested reader to~\cite{_Capria-Salamon_, salamon, _Verbitsky:HKT_, Widdows}.

\hfill

In the present note we summarise some results from the recent preprints~\cite{gra-lej-ver} and~\cite{lejmiweber}. We would like to thank the organisers of the INdAM meeting \textsl{Complex and Symplectic Geometry} for the great conference held in June 2016 in Cortona, Italy. 

%%%%%%%%%%%%%%%%%%%%%%%%%%%%%%%%%%%%%%%%%%%%%%%%%%%%%%%%%%%%

\subsection{History and examples}

While complex manifolds have been around for a long time, the study of hypercomplex manifolds only became prominent in the eighties with publications such as~\cite{boyer, salamon}. Probably the most well-known class of hypercomplex manifolds are hyperk\"{a}hler manifolds. However, the realm of hypercomplex manifolds is much broader than the one of hyperk\"{a}hler manifolds. To cite but a few hypercomplex non-hyperk\"{a}hler manifolds, note that some nilmanifolds, that is quotients of a nilpotent Lie group by a cocompact lattice, admit hypercomplex structures~\cite{BarberisDotti}. Furthermore, Dominic Joyce constructed many left-invariant hypercomplex structures on Lie groups~\cite{joyce} and similar ones have been analysed by physicists interested in string theory~\cite{SSTvP} in the context of $N=4$ supersymmetry. In more recent years, various authors constructed inhomogeneous hypercomplex structures: see for example~\cite{BoyerGalickiMann} for hypercomplex structures on Stiefel manifolds as well as \cite{battaglia, ped-poon}.

\hfill

A complete classification of compact hypercomplex manifolds of real dimension~4, called \textsl{quaternionic curves}, has been established by Charles P.~Boyer~\cite{boyer}. These are either 4-tori or K3 surfaces, both of whom are hyperk\"{a}hler, or else quaternionic Hopf surfaces~\cite{kato} which, even if non-hyperk\"{a}hler, remain locally conformally hyperk\"{a}hler. On the other hand, the situation becomes much more complicated for compact hypercomplex manifolds of real dimension~8, called \textsl{hypercomplex surfaces}. While compact complex surfaces are nowadays well understood thanks to the work of Kunihiko Kodaira~\cite{kodaira}, a similar classification for compact hypercomplex surfaces is still missing. In the sequel of this note, we will hence focus on hypercomplex manifolds of real dimension~8, the first ``unsolved dimension".

%%%%%%%%%%%%%%%%%%%%%%%%%%%%%%%%%%%%%%%%%%%%%%%%%%%%%%%%%%

\section{Cohomological properties of complex and hypercomplex manifolds}

In this Section we first briefly review some well-known cohomological aspects of complex manifolds and then show how these can be adapted to hypercomplex manifolds. For the cohomological properties of complex manifolds we refer the reader to~\cite{ang-thesis} and the references therein whilst the hypercomplex cohomologies appear in~\cite{_Capria-Salamon_, gra-lej-ver,  lejmiweber, salamon, _Verbitsky:HKT_, Widdows} to cite but a few of them. 

%%%%%%%%%%%%%%%%%%%%%%%%%%%%%%%%%%%%%%%%%%%%%%%%%%%%%%%%%%

\subsection{Cohomologies on complex manifolds}

An almost complex manifold $(X,I)$ is a smooth manifold $X$ of real dimension $2n$ together with an endomorphism of the tangent bundle $I: TX \to TX$ that satisfies $I^2=-\text{Id}_{TX}$. This \textsl{almost complex structure} $I$ can be used to decompose the bundle of complex-valued one-forms $\Omega^{1}(X) \otimes \mathbb{C}$ into the subbundle $\Omega^{1,0}_I(X)$ and the subbundle~$\Omega^{0,1}_I(X)$, with $I$ acting on the sections of $\Omega^{1,0}_I(X)$ by $i$ and on those of $\Omega^{0,1}_I(X)$ by~$-i.$ We get the following decomposition 
$$\Omega^{k}(X) \otimes \mathbb{C} = \bigoplus_{p+q=k} \Omega^{p,q}_I(X).$$ 
We denote by $\Lambda^{p,q}_I(X)$ the sections of $\Omega^{p,q}_I(X)$
and define the \textsl{Dolbeault operators}
$$\partial = \pi^{p+1,q} \circ d: \Lambda^{p,q}_I(X) \to \Lambda^{p+1,q}_I(X), \quad \bar{\partial} = \pi^{p,q+1} \circ d: \Lambda^{p,q}_I(X) \to \Lambda^{p,q+1}_I(X),$$ 
where $d$ is the exterior derivative and $\pi^{p,q}$ is the projection onto $\Lambda^{p,q}_I(X)$. Clearly, $df= \left(\partial+ \bar{\partial}\right)f$ for any function $f$. However, a priori, the same is not true for higher degree forms as explained in Figure~\ref{fig:1}:

\begin{figure}[h!]
\begin{minipage}[c]{0.35\textwidth}
\begin{equation*}
  \includegraphics{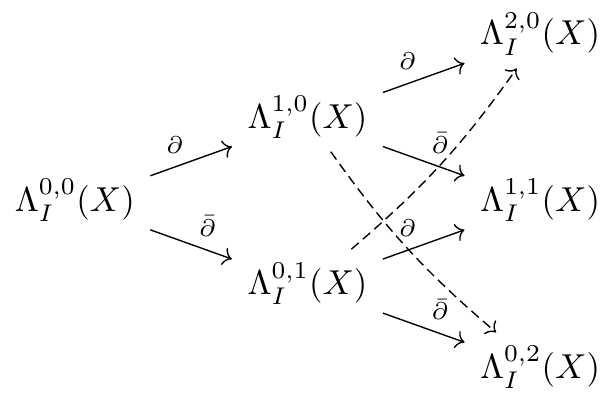}
\end{equation*}
%\begin{tikzcd}[row sep=tiny]
%&& \Lambda^{2,0}_I(X)\\				
%& \Lambda^{1,0}_I(X) \arrow[ru, "\partial"] \arrow[rd, "\bar{\partial}"] \arrow[bend right=7cm, dashed]{dddr}\\
%\Lambda^{0,0}_I(X) \arrow[rd, "\bar{\partial}"] \arrow[ru, "\partial"]  & & \Lambda^{1,1}_I(X)\\
%& \Lambda^{0,1}_I(X) \arrow[ru, "\partial"] \arrow[rd, "\bar{\partial}"] \arrow[bend right=7cm, dashed]{uuur}\\
%&& \Lambda^{0,2}_I(X)
%\end{tikzcd}
\end{minipage}\hfill
\begin{minipage}[c]{0.59\textwidth}
\caption{In general, the two dashed maps $N_I= \pi^{0,2} \circ d: \Lambda^{1,0}_I(X) \to \Lambda^{0,2}_I(X)$ and $N_I^*= \pi^{2,0} \circ d: \Lambda^{0,1}_I(X) \to \Lambda^{2,0}_I(X)$ do not need to vanish. If they do, then the almost complex structure $I$ is called integrable and $d=\partial+\bar{\partial}$ not only on functions but also on forms of higher degree.}\label{fig:1}
\end{minipage}
\end{figure}
An almost complex manifold $(X,I)$ is integrable if and only if 
\begin{equation}\label{integrability}
\partial^2 \alpha = \bar{\partial}^2 \alpha = \left(\partial \bar{\partial} + \bar{\partial} \partial\right) \alpha=0 \quad \text{ for all }\alpha \in \Lambda^{p,q}_I(X).
\end{equation} On any complex manifold $(X,I)$, there is a double complex $(\Lambda^{p,q}_I(X),\partial,\bar{\partial})$ with two anti-commuting differentials.

%%%%%%%%%%%%%%%%%%%%%%%%%%%%%%%%%%%%%%%%%%%%%%%%%%%%%%%%%%

\subsection{Cohomologies on hypercomplex manifolds}
\label{sec:3}

An \textsl{almost hypercomplex manifold} $(M,I,J,K)$ is a smooth manifold $M$ of real dimension $4n$ equipped with three almost-complex structures $I$,~$J$,~$K$ satisfying the quaternionic relations
$$I^2 = J^2 = K^2 = IJK= -\text{Id}_{TM}.$$
If all three almost-complex structures are integrable, then $(M,I,J,K)$ is called a \textsl{hypercomplex manifold}.
We would like to mimic the above characterisation of integrability~\eqref{integrability} in terms of differential operators. To this end, we will keep the decomposition of complexified differential forms with respect to the almost-complex structure~$I$. As the almost-complex structures $I$ and $J$ anticommute, we deduce that $J$ interchanges $\Lambda^{1,0}_I(M)$ with $\Lambda^{0,1}_I(M)$. This action then extends to an action $J: \Lambda^{p,q}_I(M) \to \Lambda^{q,p}_I(M)$:
$$J(\varphi)(X_1, \dots, X_p, Y_1, \dots, Y_q) = (-1)^{p+q}(\varphi)(J X_1, \dots, J X_p, J Y_1, \dots, J Y_q).$$
On any almost hypercomplex manifold, the \textsl{twisted Dolbeault operator} $\partial_J$ is defined by the commutative diagram
$$\begin{CD}
\Lambda^{p,q}_I(M)     @>\partial_J>>  \Lambda^{p+1,q}_I(M)\\
@VJVV        @AAJ^{-1}A\\
\Lambda^{q,p}_I(M)     @>\bar{\partial}>>  \Lambda^{q,p+1}_I(M)
\end{CD}$$ 
Both $\partial$ and $\partial_J$ increase the first index in the bidegree as illustrated in Figure~\ref{fig:2}:

\begin{figure}[h!]
\begin{center}
\begin{minipage}[c]{0.35\linewidth}
\begin{equation*}
  \includegraphics{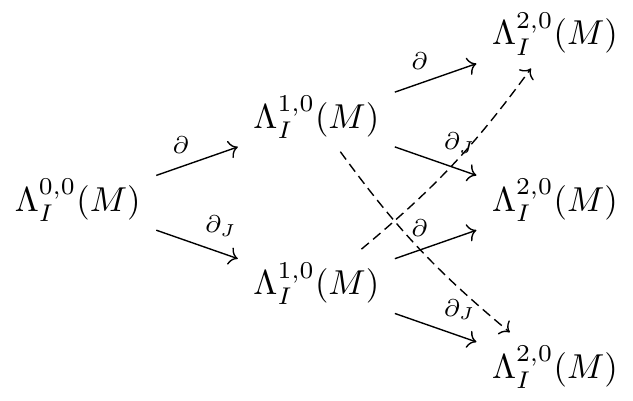}
\end{equation*}
%\begin{tikzcd}[row sep=tiny]
%&& \Lambda^{2,0}_I(M)\\				
%& \Lambda^{1,0}_I(M) \arrow[ru, "\partial"] \arrow[rd, "\partial_J"] \arrow[bend right=7cm, dashed]{dddr}\\
%\Lambda^{0,0}_I(M) \arrow[rd, "\partial_J"] \arrow[ru, "\partial"]  & & \Lambda^{2,0}_I(M)\\
%& \Lambda^{1,0}_I(M) \arrow[ru, "\partial"] \arrow[rd, "\partial_J"] \arrow[bend right=7cm, dashed]{uuur}\\
%&& \Lambda^{2,0}_I(M)
%\end{tikzcd}
\end{minipage}\hfill
\begin{minipage}[r]{0.59\linewidth}
\caption{On general almost hypercomplex manifolds, the two dashed maps $J^{-1} \circ N_I: \Lambda^{1,0}_I(M) \to \Lambda^{2,0}_I(M)$ and $N_I^* \circ J : \Lambda^{1,0}_I(M) \to \Lambda^{2,0}_I(M)$ do not need to vanish. If they do, then the almost complex structure $I$ is called integrable.}\label{fig:2}
\end{minipage}
\end{center}
\end{figure}

One checks that $\partial^2 \alpha = 0 = \partial_J^2 \alpha$ for all $\alpha \in \Lambda^{p,0}_I(M)$ if and only if the Nijenhuis tensor $N_I$ of the almost complex structure $I$ vanishes, that is if and only if the almost complex structure~$I$ is integrable. Moreover, a direct computation shows that $(\partial \partial_J + \partial_J \partial) \alpha = 0$ for all $\alpha \in \Lambda^{p,0}_I(M)$ if and only if the Nijenhuis tensor $N_J$ of the almost complex structure~$J$ vanishes. We deduce the following result~\cite{salamon, _Verbitsky:HKT_}: An almost hypercomplex manifold $(M,I,J,K)$ is integrable if and only if 
$$\partial^2 \alpha = \partial_J^2 \alpha = \left(\partial \partial_J+ \partial_J \partial \right) \alpha=0 \quad \text{ for all } \alpha \in \Lambda^{p,0}_I(M).$$ 
On any hypercomplex manifold $(M,I,J,K)$, there is always a cochain complex $(\Lambda^{p,0}_I(M),\partial,\partial_J)$ with two anti-commuting differentials. This naturally leads to a definition of cohomology groups on hypercomplex manifolds.

%%%%%%%%%%%%%%%%%%%%%%%%%%%%%%%%%%%%%%%%%%%%%%%%%%%

\subsection{Complex and Quaternionic cohomology groups}

As soon as one is facing a cochain complex with two differential operators that anticommute, one may think about defining the following cohomology groups: the Dolbeault cohomology groups, the Bott--Chern cohomology groups and the Aeppli cohomology groups. Table~\ref{tab:1} below gives precise definitions of these groups for both the double complex $(\Lambda^{p,q}_I(X),\partial, \bar{\partial})$ on a complex manifold $(X,I)$ and the single complex $(\Lambda^{p,0}_I(M),\partial,\partial_J)$ on a hypercomplex manifold $(M,I,J,K)$.

\begin{table}[h!]
\label{tab:1}  
\begin{small}
\caption{Some cohomology groups on compact complex manifolds $(X,I)$ (left) and their analogues on compact hypercomplex manifolds $(M,I,J,K)$ (right).}
\begin{tabular}{p{6.8cm}p{0.00cm}p{7.25cm}}
\hline\noalign{\smallskip}
&& \\
Complex Dolbeault cohomology groups & & Quaternionic Dolbeault cohomology groups \\
 & & \\
$H_\partial^{p,q}(X) = \frac{\{\varphi \in \Lambda^{p,q}_I(X)~|~\partial \varphi=0\}}{\partial \Lambda^{p-1,q}_I(X)} = \frac{\text{Ker }\partial}{\text{Im }\partial}$ & & $H_\partial^{p,0}(M) = \frac{\{\varphi \in \Lambda^{p,0}_I(M)~|~\partial \varphi=0\}}{\partial \Lambda^{p-1,0}_I(M)} = \frac{\text{Ker }\partial}{\text{Im }\partial}$ \\
 & & \\
 $H_{\bar{\partial}}^{p,q}(X) = \frac{\{\varphi \in \Lambda^{p,q}_I(X)~|~\bar{\partial} \varphi=0\}}{\bar{\partial} \Lambda^{p,q-1}_I(X)} = \frac{\text{Ker }\bar{\partial}}{\text{Im }\bar{\partial}}$ & & $H_{\partial_J}^{p,0}(M) = \frac{\{\varphi \in \Lambda^{p,0}_I(M)~|~\partial_J \varphi=0\}}{\partial_J \Lambda^{p-1,0}_I(M)} = \frac{\text{Ker }\partial_J}{\text{Im }\partial_J}$ \\
& & \\
\noalign{\smallskip}\hline\noalign{\smallskip}
&& \\
Complex Bott--Chern cohomology groups & & Quaternionic Bott--Chern cohomology groups \\
& & \\
$H^{p,q}_{BC}(X) = \frac{\{\varphi \in \Lambda^{p,q}_I(X)~|~\partial \varphi = 0 = \bar{\partial} \varphi \}}{\partial \bar{\partial} \Lambda^{p-1,q-1}_I(X)} = \frac{\text{Ker }\partial \cap \text{Ker }\bar{\partial}}{\text{Im }\partial \bar{\partial}}$ & & $H^{p,0}_{BC}(M) = \frac{\{\varphi \in \Lambda^{p,0}_I(M)~|~\partial \varphi = 0 = \partial_J \varphi \}}{\partial \partial_J \Lambda^{p-2,0}_I(M)} = \frac{\text{Ker }\partial \cap \text{Ker }\partial_J}{\text{Im }\partial \partial_J}$ \\
 & & \\
\noalign{\smallskip}\hline\noalign{\smallskip}
&& \\
Complex Aeppli cohomology groups & & Quaternionic Aeppli cohomology groups \\
& & \\
$H^{p,q}_{AE}(X) = \frac{\{\varphi \in \Lambda^{p,q}_I(X)~|~\partial \bar{\partial} \varphi = 0\}}{\partial \Lambda^{p-1,q}_I(X) + \bar{\partial} \Lambda^{p,q-1}_I(X)} =  \frac{\text{Ker }\partial \bar{\partial}}{\text{Im }\partial + \text{Im }\bar{\partial}}$ & & $H^{p,0}_{AE}(M) = \frac{\{\varphi \in \Lambda^{p,0}_I(M)~|~\partial \varphi = 0 = \partial_J \varphi \}}{\partial \partial_J \Lambda^{p-2,0}_I(M)} = \frac{\text{Ker }\partial \cap \text{Ker }\partial_J}{\text{Im }\partial \partial_J}$ \\
 & & \\
 \noalign{\smallskip}\hline\noalign{\smallskip}
\end{tabular}
\end{small}
\end{table}

On compact hypercomplex manifolds, the groups $H^{p,0}_\partial(M)$, $H^{p,0}_{\partial_J}(M)$, $H^{p,0}_{BC}(M)$ and $H^{p,0}_{AE}(M)$ are finite-dimensional complex vector spaces~\cite{gra-lej-ver}. 

%%%%%%%%%%%%%%%%%%%%%%%%%%%%%%%%%%%%%%%%%%%%%%%%%%%

\subsection{Conjugation symmetry}

On a complex manifold $(X,I)$, conjugation defines a map 
$$\Lambda^{p,q}_I(X) \to \Lambda^{q,p}_I(X) : \alpha \mapsto \bar{\alpha}.$$ 
As this map passes to cohomology, we deduce that $H^{p,q}_{\partial}(X) \cong H^{q,p}_{\bar{\partial}}(X)$. Furthermore, this also implies that 
$$H^{p,q}_{BC}(X) \cong H^{q,p}_{BC}(X) \quad \text{ and } \quad H^{p,q}_{AE}(X) \cong H^{q,p}_{AE}(X).$$
On a hypercomplex manifold $(M,I,J,K)$, conjugation followed by the action of $J$ similarly defines a map 
$$\Lambda^{p,0}_I(M) \to \Lambda^{p,0}_I(M) : \alpha \mapsto J(\bar{\alpha}).$$
Once more, this map descends to cohomology and leads to the isomorphism
$$H^{p,0}_{\partial}(M) \cong H^{p,0}_{\partial_J}(M)$$
but we do not get any isomorphisms for $H^{p,0}_{BC}(M)$ or $H^{p,0}_{AE}(M)$.

%%%%%%%%%%%%%%%%%%%%%%%%%%%%%%%%%%%%%%%%%%%%%%%%%%%

\subsection{The $\partial\partial_J$-Lemma}

On a compact complex manifold $(X,I)$ and on a compact hypercomplex manifold $(M,I,J,K)$, the identity map induces the following maps: 

\begin{center}
\begin{equation*}
  \includegraphics{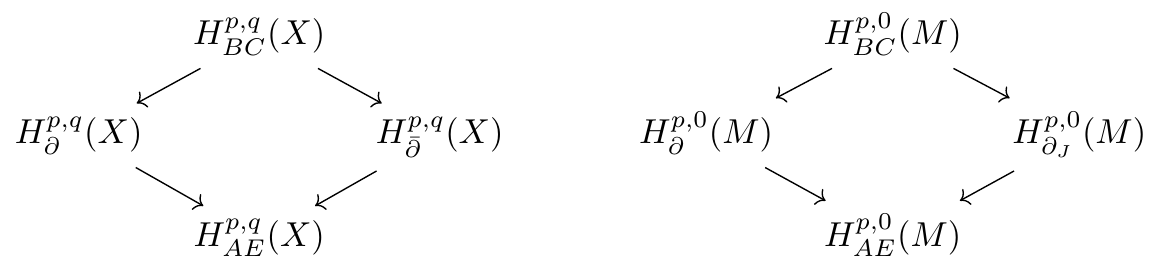}
\end{equation*}
%\begin{tikzcd}[column sep=tiny, row sep=small]
%& H^{p,q}_{BC}(X) \arrow[ld] \arrow[rd] & & \quad & & H^{p,0}_{BC}(M) \arrow[ld] \arrow[rd] &  \\		
%H^{p,q}_{\partial}(X) \arrow[rd] & & H^{p,q}_{\bar{\partial}}(X) \arrow[ld] & \quad & H^{p,0}_{\partial}(M) \arrow[rd] & & H^{p,0}_{\partial_J}(M) \arrow[ld] \\
%& H^{p,q}_{AE}(X) & & \quad & & H^{p,0}_{AE}(M) & 
%\end{tikzcd}
\end{center}

In general, these maps have no reason to be either injective or surjective. We say that the $\partial\bar{\partial}$-Lemma holds if the map $H^{p,q}_{BC}(X) \to H^{p,q}_{\bar{\partial}}(X)$ is injective and similarly that the $\partial \partial_J$-Lemma is satisfied if the map $H_{BC}^{p,0}(M) \to H_{\partial_J}^{p,0}(M)$ is injective. In other words, the $\partial\bar{\partial}$-Lemma holds if every $\partial$-closed $\bar{\partial}$-exact $(p,q)$-form is $\partial \bar{\partial}$-exact while the $\partial \partial_J$-Lemma holds if every $\partial$-closed, $\partial_J$-exact $(p,0)$-form is $\partial \partial_J$-exact. As it turns out, this actually implies that all of the maps in the above diagram become isomorphisms~\cite{DGMS}. 

%%%%%%%%%%%%%%%%%%%%%%%%%%%%%%%%%%%%%%%%%%%%%%%%%%%

\subsection{A quaternionic Fr\"{o}licher-type inequality}

We deduce that, on a compact complex manifold $(X,I)$, the Bott--Chern and Aeppli cohomology groups may differ from the Dolbeault and deRham cohomology groups (if the $\partial \bar{\partial}$-Lemma does not hold). The following result by Angella--Tomassini quantifies this difference:

\begin{theorem}\label{theorem1}\cite{_Angella-Tomassini_}
Let $(X,I)$ be a compact complex manifold of real dimension~$2n$. Then
\begin{equation}\label{ineqComplex}
\sum_{p+q=k} \left( \dim H^{p,q}_{BC}(X) + \dim H^{p,q}_{AE}(X) \right) \geqslant 2 \dim H^k_{dR}(X) 
\end{equation}
for any $0\leqslant k\leqslant n$ where
$$H^{k}_{dR}(X) = \frac{\{\varphi \in \Lambda^{k}(X)~|~d\varphi=0\}}{d \Lambda^{k}(X)} = \frac{\text{Ker }d}{\text{Im }d}$$
denotes deRham cohomology. Moreover, the $\partial \bar{\partial}$-Lemma holds if and only if we have equality for all $0 \leqslant k \leqslant n$.
\end{theorem}

A similar result can be established for quaternionic cohomologies on compact hypercomplex manifolds:

\begin{theorem}\label{theorem2}\cite{lejmiweber}
Let $(M,I,J,K)$ be a compact hypercomplex manifold of real dimension~$4n$. Then
\begin{equation}\label{ineqHypercomplex}
\dim H_{BC}^{p,0}(M) + \dim H_{AE}^{p,0}(M) \geqslant 2 \dim E_2^{p,0}(M)
\end{equation}
for any $0 \leqslant p \leqslant 2n$ where the space $E_2^{p,0}(M)$ is defined by
$$E_2^{p,0}(M) = \frac{ \{\varphi \in \Lambda^{p,0}_I(M)~|~\partial \varphi = 0 \text{ and } \partial_J \varphi +\partial \alpha_1 =0 \}}{\{\varphi \in \Lambda^{p,0}_I(M)~|~\varphi=\partial \beta_1 + \partial_J \beta_2 \text{ and } \partial \beta_2 =0\}}.$$
Moreover, the $\partial \partial_J$-Lemma holds if and only if we have equality for all $0\leqslant p \leqslant 2n$.
\end{theorem}

While these results look very similar, the conclusions we draw differ. More precisely, recall that the Betti numbers appearing in the right-hand-side of~\eqref{ineqComplex} are topological invariants. As the dimensions of the cohomology groups are upper semi-continuous, Angella and Tomassini deduce from Theorem~\ref{theorem1} that, on compact complex manifolds, the $\partial \bar{\partial}$-Lemma is stable by small complex deformations~\cite{_Angella-Tomassini_, voisin, wu}. However, the same reasoning fails on compact hypercomplex manifolds, because the term $\dim E_2^{p,0}(M)$ appearing in the right-hand-side of~\eqref{ineqHypercomplex} in Theorem~\ref{theorem2} has no reason to be a topological invariant. Indeed, it can be shown that the $\partial \partial_J$-Lemma is not stable by small hypercomplex deformations as illustrated in the Example in Section~\ref{SectionExample}.\\

Finally, Theorems~\ref{theorem1} and~\ref{theorem2} also allow us to quantify how far away a complex manifold is from being ``cohomologically K\"{a}hler" and similarly how far away a hypercomplex manifold is from being ``cohomologically HKT" (see Section~\ref{SectionMetrics}). Define the \textsl{non-K\"{a}hler-ness degrees}~\cite{ang-dlo-tom} on complex manifolds
$$\Delta^k(X) = \sum_{p+q=k} \left(\dim H^{p,q}_{BC}(X) + \dim H^{p,q}_{AE}(X)\right) - 2 \dim H^{k}_{dR}(X)$$
and the \textsl{non-HKT-ness degrees}~\cite{lejmiweber} on hypercomplex manifolds
$$\Delta^p(M) = \dim H^{p,0}_{BC}(M) + \dim H^{p,0}_{AE}(M) - 2 \dim E^{p,0}_2(M).$$

%%%%%%%%%%%%%%%%%%%%%%%%%%%%%%%%%%%%%%%%%%%%%%%%%%%

\section{Metric structures}\label{SectionMetrics}

Every complex manifold $(X,I)$ admits a Hermitian metric, that is~a Riemannian metric~$g$ such that
$$g(\cdot,\cdot)= g(I\cdot,I\cdot).$$
We can build out of this the Hermitian form $\omega(\cdot,\cdot)=g(I\cdot,\cdot)$ and various special metrics can be characterised via conditions on $\omega$. Similarly, any hypercomplex manifold $(M,I,J,K)$ admits a quaternionic Hermitian metric, that is a Riemannian metric $g$ which satisfies
$$g(\cdot,\cdot) = g(I\cdot,I\cdot) = g(J\cdot, J\cdot) =g(K\cdot,K\cdot).$$
This leads to three (not necessarily closed) differential forms $\omega_I(\cdot,\cdot)=g(I\cdot,\cdot)$, $\omega_J(\cdot,\cdot)= g(J\cdot,\cdot)$ and $\omega_K(\cdot,\cdot)= g(K\cdot,\cdot)$
that can be assembled to build the fundamental form
$$\Omega= \omega_J + \sqrt{-1} \omega_K$$
which is of type $(2,0)$ with respect to the complex structure $I$. Once more, various special metrics can be characterised by imposing conditions on the form $\Omega$. If, for instance, the form~$\Omega$ is $d$-closed then $(M,I,J,K,\Omega)$ is called a \textsl{hyperk\"{a}hler} manifold whereas if $\Omega$ is $\partial$-closed, then $(M,I,J,K,\Omega)$ is called \textsl{hyperk\"{a}hler with torsion}, or \textsl{HKT} for short (see~\cite{gra-poon} for a nice introduction). Table~\ref{tab:2} summarises some special metrics on hypercomplex manifolds together with their associated conditions on $\Omega$ as well as their complex counterparts. We point out that HKT metrics, just as K\"{a}hler metrics in the complex setup, admit a local potential~\cite{BanosSwann}. 

\begin{table}[h!]
\caption{Correspondence between metric structures on complex and hypercomplex manifolds}
\label{tab:2} 
\begin{tabular}{p{3.2cm}p{2.25cm}p{5.2cm}p{2.7cm}}
\hline\noalign{\smallskip}
Complex & Condition & Hypercomplex & Condition  \\
&&& \\
\noalign{\smallskip}\hline\noalign{\smallskip}
Gauduchon & $\partial\bar{\partial}\omega^{n-1}=0$  & Quaternionic Gauduchon & $\partial\partial_J\Omega^{n-1}=0$\\
Strongly Gauduchon & $\partial \omega^{n-1} \in \text{Im }{\bar{\partial}}$ & Quaternionic strongly Gauduchon & $\partial \Omega^{n-1} \in \text{Im }\partial_J$\\
Balanced & $d\omega^{n-1}=0$ & Quaternionic balanced & $\partial \Omega^{n-1}=0$ \\
K\"{a}hler & $d\omega=0$ & Hyperk\"{a}hler with torsion (HKT) & $\partial \Omega=0$ \\
 &  & Hyperk\"{a}hler & $d\Omega=0$\\
 &&& \\
\noalign{\smallskip}\hline\noalign{\smallskip}
\end{tabular}
\end{table}

A first important difference between complex and hypercomplex manifolds is the existence of a preferred metric. Indeed, a complex manifold always admits a Gauduchon metric and this metric is unique in its conformal class up to a constant. On the other hand, to recover existence of a quaternionic Gauduchon metric on hypercomplex manifolds, we will impose an additional holonomy constraint as described in the next Section.

%%%%%%%%%%%%%%%%%%%%%%%%%%%%%%%%%%%%%%%%%%%%%%%%%%%

\section{$SL(n,\mathbb{H})$-manifolds}
\label{sec:4}

There is a particular class of hypercomplex manifolds, called $SL(n,\mathbb{H})$-manifolds, that shares more properties of complex manifolds than general hypercomplex manifolds do. The key reason for this is that the canonical bundle of an $SL(n,\mathbb{H})$-manifold is holomorphically trivial and this leads to a version of Hodge theory when HKT~\cite{_Verbitsky:HKT_} and to a version of Serre duality on the bundle~$\Omega^{*,0}_I(M)$.

%%%%%%%%%%%%%%%%%%%%%%%%%%%%%%%%%%%%%%%%%%%%%%%%%%%

\subsection{The Obata connection}

\begin{figure}[h!]
\begin{minipage}[c]{0.45\textwidth}
\caption{This Figure shows the possible holonomy groups of compact hypercomplex manifolds in real dimension~8. Left-invariant structures on Lie groups are conjectured to have holonomy equal to $GL(2,\mathbb{H})$ just as it has been proven for $SU(3)$.}\label{fig:3}
\end{minipage} \hfill
\begin{minipage}[c]{0.5\textwidth}
\begin{equation*}
  \includegraphics{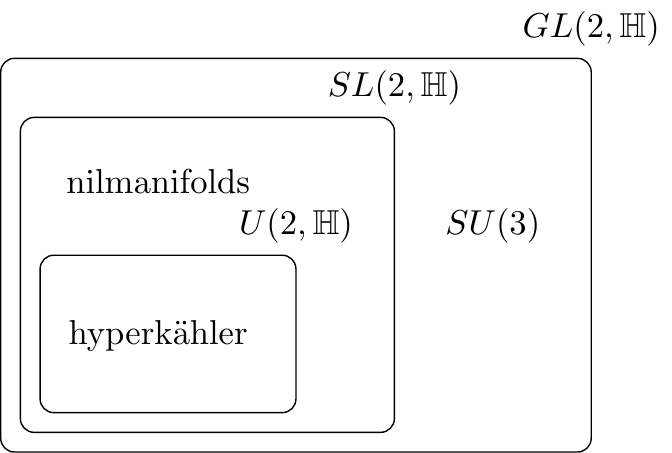}
\end{equation*}
%\begin{tikzpicture}[rounded corners, fill=gray]
%\draw 
%      (2,0) node [text=black,right,above] {$SU(3)$}
%      (-1.4,0.5) node [text=black,right,above] {nilmanifolds}
%      (-1.4,-1.1) node [text=black,right,above] {hyperk\"{a}hler}
%      (-2.8,-1.8) rectangle (1,1.4) (1,1.4)  node [text=black,right,above] {$SL(2,\mathbb{H})$}
%      (-2.6,-1.6) rectangle (0,0) (0,0)  node [text=black,right, above] {$U(2,\mathbb{H})$}
%      (-3,-2) rectangle (3,2) node [text=black,right,above] {$GL(2,\mathbb{H})$};  
%\end{tikzpicture}
\end{minipage}
\end{figure}

Another important difference between complex and hypercomplex geometry is the existence of a special connection. A complex manifold generally admits infinitely many torsion-free connections which preserve the complex structure~\cite{JoyceBook}. 
On the other hand, any hypercomplex manifold admits a unique torsion-free connection $\nabla$ such that
$$\nabla I = \nabla J = \nabla K = 0.$$
This connection is called the Obata connection~\cite{obata}. In general, the Obata connection does not preserve the metric, except when the manifold is hyperk\"{a}hler. Given any torsion-free affine connection, the holonomy group introduced by \'{E}lie Cartan measures the failure of the parallel translation associated to a connection to be holonomic. Merkulov and Schwachh\"{o}fer classified the groups which can possibly arise as irreducible holonomy groups of torsion-free connections~\cite{MerkulovSchwachhofer}. As illustrated in Figure~\ref{fig:3}, in the case of the Obata connection, there are three possible choices: $GL(n,\mathbb{H})$, $SL(n,\mathbb{H})$ and $U(n,\mathbb{H})$.
Indeed, as the Obata connection preserves all three complex structures, its holonomy is necessarily contained in the quaternionic general linear group~$GL(n,\mathbb{H})$. It turns out that, for all nilmanifolds, the holonomy is contained in the commutator subgroup $SL(n,\mathbb{H})$~\cite{_BDV:nilmanifolds_}. Finally, a hyperk\"{a}hler manifold is characterised by the fact that the holonomy of the Obata connection is equal to the compact symplectic group $Sp(n)=U(n,\mathbb{H})$, that is the hyperunitary group. For the homogeneous hypercomplex structure on $SU(3)$ constructed by Joyce, the holonomy is equal to~$GL(2,\mathbb{H})$~\cite{_Soldatenkov:SU(2)_}.

%%%%%%%%%%%%%%%%%%%%%%%%%%%%%%%%%%%%%%%%%%%%%%%%%%%

\subsection{Hodge theory}

One important aspect of $SL(n,\mathbb{H})$-manifolds is that, if the metric is HKT, then it is possible to establish a version of Hodge theory~\cite{_Verbitsky:HKT_}. Indeed, any $SL(n,\mathbb{H})$-manifold has holomorphically trivial canonical bundle. The nowhere degenerate real holomorphic section $\Phi$ (that is a nowhere degenerate section $\Phi$ such that $J\bar{\Phi}=\Phi$ and $\partial \bar{\Phi} = 0$) which trivialises $\Omega^{2n,0}_I(M)$ may then be used to define a Hodge star operator on a $SL(n,\mathbb{H})$-manifold $(M,I,J,K,\Omega)$
$$\star_\Phi : \Lambda^{p,0}_I(M) \to \Lambda^{2n-p,0}_I(M)$$
via
$$\alpha \wedge (\star_\Phi \beta) \wedge \bar{\Phi} = h(\alpha, \beta) \frac{\Omega^n \wedge \bar{\Phi}}{n!},$$
where $h$ is the $\mathbb{C}$-bilinear extension with respect to $I$ of the quaternionic Hermitian metric $g$ associated to $\Omega$. On compact manifolds, this leads to the adjoints 
$$\partial^{*_\Phi}= - \star_\Phi \partial \star_\Phi \quad \text{ and } \quad \partial_J^{\star_\Phi} = - \star_{\Phi} \partial_J \star_{\Phi}$$
and thus to the Laplacians
$$\Delta_\partial = \partial \partial^{*_\Phi} + \partial^{*_\Phi} \partial \quad \text{ and } \quad \Delta_{\partial_J} = \partial_J \partial_J^{*_\Phi} + \partial_J^{*_\Phi} \partial_J.$$
On $SL(2,\mathbb{H})$-manifolds, the Hodge $\star_{\Phi}$ acts as an involution on $\Lambda^{2,0}_I(M)$ and hence we may decompose $(2,0)$-forms into $\star_\Phi$-self-dual ones and $\star_\Phi$-anti-self-dual ones. As $\star_{\Phi}$ commutes with~$\Delta_\partial$, this splitting descends to cohomology. We conclude that, on a compact $SL(2,\mathbb{H})$-manifold, the space~$H^{2,0}_{\partial}(M)$ can be decomposed as a direct sum of $\partial$-closed $\star_\Phi$-self-dual and $\partial$-closed $\star_\Phi$-anti-self-dual forms. 

%%%%%%%%%%%%%%%%%%%%%%%%%%%%%%%%%%%%%%%%%%%%%%%%%%%

\subsection{Serre duality and $SL(n,\mathbb{H})$-symmetry}

\begin{figure}[h!]
\begin{center}
\begin{equation*}
  \includegraphics{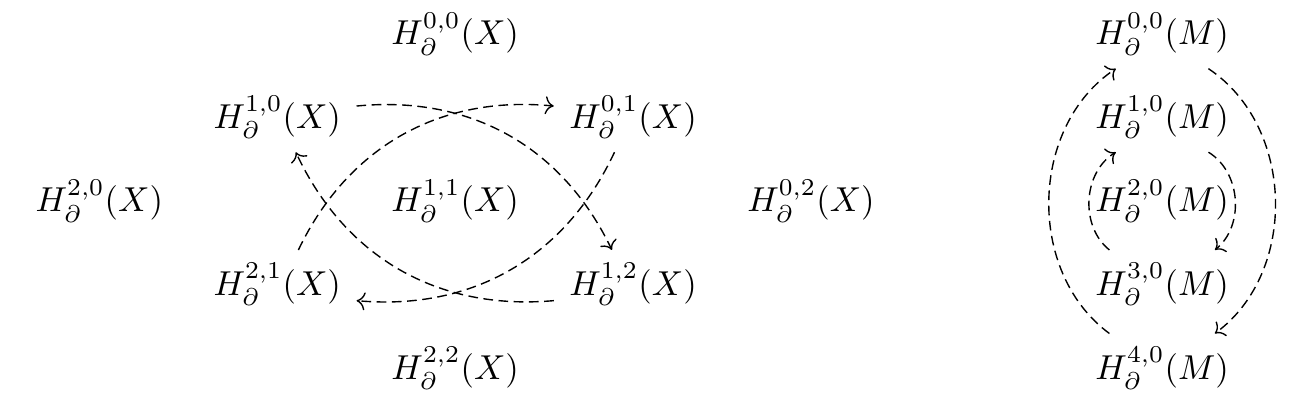}
\end{equation*}
%\begin{tikzcd}[column sep=tiny, row sep=tiny]
%& & & H^{0,0}_{\partial}(X) & \quad & \quad & \quad & \quad & H^{0,0}_{\partial}(M) \arrow[dashed,bend left=55]{dddd}\\		 
%& & H^{1,0}_{\partial}(X) \arrow[dashed,bend left=35]{rrdd}& \quad &H^{0,1}_{\partial}(X) \arrow[dashed,bend left=35]{ddll} & \quad & \quad & \quad & H^{1,0}_{\partial}(M) \arrow[dashed,bend left=55]{dd}\\
%& H^{2,0}_{\partial}(X) & \quad & H^{1,1}_{\partial}(X) & \quad & H^{0,2}_{\partial}(X) & \quad & \quad & H^{2,0}_{\partial}(M)\\
%& \quad & H^{2,1}_{\partial}(X) \arrow[dashed,bend left=35]{uurr}& \quad & H^{1,2}_{\partial}(X) \arrow[dashed,bend left=35]{lluu}& \quad & \quad & \quad & H^{3,0}_{\partial}(M) \arrow[dashed,bend left=55]{uu}\\
%& \quad & \quad & H^{2,2}_{\partial}(X) & \quad & \quad & \quad & \quad & H^{4,0}_{\partial}(M) \arrow[dashed,bend left=55]{uuuu}
%\end{tikzcd}
\end{center}
\caption{Serre duality on compact complex surfaces (left) and $SL(2,\mathbb{H})$-symmetry on compact hypercomplex surfaces (right).}\label{fig:4}
\end{figure}

Besides the conjugation symmetry, compact complex manifolds also satisfy Serre duality coming from the pairing on $H^{p,q}_\partial(X) \times H^{n-p,n-q}_{\partial}(X)$ given by
$$([\alpha],[\beta]) \mapsto \int_X \alpha \wedge \beta.$$
On compact $SL(n,\mathbb{H})$ manifolds, an analogue of this exists and can be formulated as follows (see also Figure~\ref{fig:4}). Consider the pairing $H^{p,0}_{\partial}(M) \times H^{2n-p,0}_{\partial}(M)$ given by
$$([\alpha],[\beta]) \mapsto \int_M \alpha \wedge \beta \wedge \bar{\Phi}.$$
Note that, for this to be well-defined we really need $\partial \bar{\Phi}=0$.

Furthermore, Serre duality and $SL(n,\mathbb{H})$-symmetry also provide links between Bott--Chern and Aeppli cohomologies. Indeed, using the above pairings, it can be shown that Serre duality on compact complex manifolds of real dimension $2n$ implies that $H^{p,q}_{BC}(X) \cong H^{n-p,n-q}_{AE}(X)$ and similarly $SL(n,\mathbb{H})$-symmetry on compact $SL(n,\mathbb{H})$-manifolds implies that $H^{p,0}_{BC}(M) \cong H^{2n-p,0}_{AE}(M)$.

%%%%%%%%%%%%%%%%%%%%%%%%%%%%%%%%%%%%%%%%%%%%%%%%%%%

\subsection{$SL(2,\mathbb{H})$-manifolds}

We saw that compact $SL(2,\mathbb{H})$-manifolds share many properties of compact complex surfaces, most notably a version of Hodge theory when it is HKT and similar symmetries. Hence it should not surprise that many results valid on compact complex surfaces can be adapted to results on $SL(2,\mathbb{H})$-manifolds. To illustrate this link, we display in Table~\ref{tab:3} some results which show that HKT metrics play a similar role on $SL(2,\mathbb{H})$-manifolds than K\"{a}hler metrics do on complex surfaces.

\begin{table}[h!]
\caption{Results valid on compact complex surfaces (left) and the corresponding results on compact $SL(2,\mathbb{H})$-manifolds (right).}
\label{tab:3}  
\begin{tabular}{p{5.8cm}p{1.6cm}p{6.2cm}}
\hline\noalign{\smallskip}
Compact complex surfaces & & Compact $SL(2,\mathbb{H})$-manifolds   \\
 & & \\
\noalign{\smallskip}\hline\noalign{\smallskip}
K\"{a}hler if and only if $\dim H^1_{dR}(X)$ even~\cite{buchdahl, lamari, miyaoka, siu} & & HKT if and only if $\dim H^{1,0}_{\partial}(M)$ even~\cite{gra-lej-ver} \\
 & & \\
\noalign{\smallskip}\hline\noalign{\smallskip}
K\"{a}hler if and only if strongly Gauduchon~\cite{pop} & & HKT if and only if quaternionic strongly Gauduchon~\cite{lejmiweber} \\
 & & \\
\noalign{\smallskip}\hline\noalign{\smallskip}
K\"{a}hler if and only if the second non-K\"{a}hler-ness degree vanishes~\cite{ang-dlo-tom, lub-tel, teleman} & & HKT if and only if the second non-HKT-ness degree vanishes~\cite{lejmiweber} \\
 & & \\
 \noalign{\smallskip}\hline\noalign{\smallskip}
\end{tabular}
\end{table}

%%%%%%%%%%%%%%%%%%%%%%%%%%%%%%%%%%%%%%%%%%%%%%%%%%%

\subsection{Computations}\label{SectionExample}

On a compact hypercomplex nilmanifold $(M,I,J,K)$ of real dimension~8, if one assumes that the Dolbeault cohomology $H_{\bar{\partial}}^{p,q}(X)$ with respect to $I$ can be computed
using left-invariants forms then the quaternionic Dolbeault cohomology groups~$H^{p,0}_{\partial}(M)$ and~$H^{p,0}_{\partial_J}(M)$, the quaternionic Bott--Chern cohomology groups $H^{p,0}_{BC}(M)$ as well as the quaternionic Aeppli cohomology groups $H^{p,0}_{AE}(M)$ can be computed using only left-invariant forms~\cite{lejmiweber}. Hence we may explicitly calculate these cohomologies for the following example based upon the central extension of the quaternionic Lie algebra $\mathbb{R} \times H_7$. We consider a path of hypercomplex structures as done in~\cite{fin-gra, gra-lej-ver, lejmiweber}. We end up with an $SL(2,\mathbb{H})$-manifold carrying a family $t\in (0,1)$ of hypercomplex structures which is HKT for $t=\frac{1}{2}$ but not HKT for all other values of $t$. The structure equations of the Lie algebra are:
$$\left\{
\begin{array}{l}
de^1 = de^2= de^3 =de^4 = de^5 =0, \\
de^6= e^1 \wedge e^2 + e^3 \wedge e^4, \\
de^7= e^1 \wedge e^3 + e^4 \wedge e^2, \\
de^8= e^1 \wedge e^4 + e^2 \wedge e^3.
\end{array}
\right.$$
Consider the family of hypercomplex structures $(I_t,J_t,K_t)$ parametrised by $t \in (0,1)$:
\begin{equation*}
\begin{array}{ccccccc}
I_t e^1= \frac{t-1}{t} e^2, & \quad & I_t e^3 = e^4, & \quad & I_t e^5 = \frac{1}{t} e^6, & \quad & I_t e^7 = e^8, \\
J_t e^1= \frac{t-1}{t} e^3, & \quad & J_t e^2 = - e^4, & \quad & J_t e^5 = \frac{1}{t} e^7, & \quad & J_t e^6 = - e^8. 
\end{array}
\end{equation*}
A basis of left-invariant $(1,0)$-forms is given by:
\begin{equation*}
\begin{array}{ccccccc}
\varphi^1 = e^1 - i \frac{t-1}{t} e^2, & \quad & \varphi^2 = e^3- i e^4, & \quad & \varphi^3 = e^5 - i \frac{1}{t} e^6, & \quad & \varphi^4 = e^7 - i e^8.
\end{array}
\end{equation*}
The structure equations become:
\begin{equation*}
\begin{array}{ccccccc}
d\varphi^1 = 0, & \quad & d\varphi^2 = 0, & \quad & d\varphi^3 =  \frac{1}{2(1-t)} \varphi^{1\bar{1}} - \frac{1}{2t} \varphi^{2\bar{2}}, & \quad & d\varphi^4 = \frac{2t-1}{2t-2} \varphi^{12} - \frac{1}{2t-2} \varphi^{\bar{1}2}.
\end{array}
\end{equation*}
If $t =\frac{1}{2}$, then $d\varphi^{i} \subseteq \Lambda^{1,1}_I(M)$ and the complex structure is abelian whereas otherwise it is not.
In terms of the differentials $\partial$ and $\partial_J$, the structure equations can be written as:
\begin{equation*}
\begin{array}{rcrclclcc}
\partial \varphi^1 = 0, & \quad & \partial \varphi^2 = 0, & \quad & \partial \varphi^3 = 0, & \quad & \partial \varphi^4 = \frac{2t-1}{2(t-1)} \varphi^{12}, \\
\partial_J \varphi^1 = 0, & \quad & \partial_J \varphi^2 = 0, & \quad & \partial_J \varphi^3 = - \frac{2t-1}{2(t-1)}\varphi^{12}, & \quad & \partial_J \varphi^4 = 0.
\end{array}
\end{equation*}
We conclude: if $t\neq \frac{1}{2}$, then we get Table~\ref{tab:4}:

\begin{table}[h!]
\caption{Dimensions of the quaternionic cohomology groups when $t=\frac{1}{2}$.}
\label{tab:4} 
\begin{center}
\begin{tabular}{p{1.5cm}p{1cm}p{1cm}p{1cm}p{1cm}}
\hline\noalign{\smallskip}
$(p,0)$ & $h_\partial^{p,0}$ & $h^{p,0}_{\partial_J}$ & $h^{p,0}_{BC}$ & $h^{p,0}_{AE}$ \\
\noalign{\smallskip}\hline\noalign{\smallskip}
$(1,0)$ & $3$ & $3$ & $2$ & $4$ \\
$(2,0)$ & $4$ & $4$ & $5$ & $5$ \\
$(3,0)$ & $3$ & $3$ & $4$ & $2$ \\
\noalign{\smallskip}\hline\noalign{\smallskip}
\end{tabular}
\end{center}
\end{table}
\noindent whereas if $t=\frac{1}{2}$ then both $\partial \varphi^4=0$ and $\partial_J \varphi^3=0$ which leads to Table~\ref{tab:5}:

\begin{table}[h!]
\caption{Dimensions of the quaternionic cohomology groups when $t \neq \frac{1}{2}$.}
\label{tab:5}
\begin{center}
\begin{tabular}{p{1.5cm}p{1cm}p{1cm}p{1cm}p{1cm}}
\hline\noalign{\smallskip}
$(p,0)$ & $h_\partial^{p,0}$ & $h^{p,0}_{\partial_J}$ & $h^{p,0}_{BC}$ & $h^{p,0}_{AE}$ \\
\noalign{\smallskip}\hline\noalign{\smallskip}
$(1,0)$ & $4$ & $4$ & $4$ & $4$ \\
$(2,0)$ & $6$ & $6$ & $6$ & $6$ \\
$(3,0)$ & $4$ & $4$ & $4$ & $4$ \\
\noalign{\smallskip}\hline\noalign{\smallskip}
\end{tabular}
\end{center}
\end{table}

We deduce that, just as the HKT property, the $\partial \partial_J$-Lemma is not stable by small hypercomplex deformations~\cite{fin-gra}. This differs from the complex setup where the $\partial \bar{\partial}$-Lemma is stable by small complex deformations~\cite{_Angella-Tomassini_, voisin, wu}.

%%%%%%%%%%%%%%%%%%%%%%%%%%%%%%%%%%%%%%%%%%%%%%%%%%%

\end{document}